
%
%
\def\rtimes{\#}
\def\z{{_{(0)}}}
\def\o{{_{(1)}}}
\def\t{{_{(2)}}}
\def\th{{_{(3)}}}

\def\id{{\rm id}}
\def\ker{{\rm ker}}

\def\d{{\rm d}}
\def\ad{{\rm Ad_{R}}}
\def\S{{\rm S}}
\def\CF{{\cal F}}

\def\CN{{\cal N}}
\def\CQ{{\cal Q}}
\def\CH{{\cal H}}
\def\C{{\bf C}}
\def\P{{\bf P}}
\def\du{{\rm d_ U}}

\def\tens{{\otimes}}

\def\endproof{{\ $\Box$}\bigskip }
\def\eps{{\epsilon}}
\def\tauo{{\tau^{(1)}}}\def\taut{{\tau^{(2)}}}

\def\omp{{\Omega^{1}(P) }}
\def\omup{{\Omega^{1}P}}


\frenchspacing

\parindent15pt

\abovedisplayskip4pt plus2pt
\belowdisplayskip4pt plus2pt
\abovedisplayshortskip2pt plus2pt
\belowdisplayshortskip2pt plus2pt

\font\twbf=cmbx10 at12pt
 at12pt
 at12pt

\font\sc=cmcsc10

\font\ninerm=cmr9
\font\nineit=cmti9
\font\ninesy=cmsy9
\font\ninei=cmmi9
\font\ninebf=cmbx9

\font\sevenrm=cmr7

\font\seveni=cmmi7
\font\sevensy=cmsy7

\font\fivenrm=cmr5
\font\fiveni=cmmi5
\font\fivensy=cmsy5

\def\nine{%
\textfont0=\ninerm \scriptfont0=\sevenrm \scriptscriptfont0=\fivenrm
\textfont1=\ninei \scriptfont1=\seveni \scriptscriptfont1=\fiveni
\textfont2=\ninesy \scriptfont2=\sevensy \scriptscriptfont2=\fivensy
\textfont3=\tenex \scriptfont3=\tenex \scriptscriptfont3=\tenex
\def\rm{\fam0\ninerm}%
\textfont\itfam=\nineit
\def\it{\fam\itfam\nineit}%
\textfont\bffam=\ninebf
\def\bf{\fam\bffam\ninebf}%
\normalbaselineskip=11pt
\setbox\strutbox=\hbox{\vrule height8pt depth3pt width0pt}%
\normalbaselines\rm}

\hsize30cc
\vsize44cc
\nopagenumbers

\def\luz#1{\luzno#1?}
\def\luzno#1{\ifx#1?\let\next=\relax\yyy
\else \let\next=\luzno#1\xxx\fi\next}
\def\sp#1{\def\xxx{\kern1.7pt}\def\yyy{\kern-1.7pt}\luz{#1}}
\def\spa#1{\def\xxx{\kern1pt}\def\yyy{\kern-1pt}\luz{#1}}

\newcount\beg
\newbox\aabox
\newbox\atbox
\newbox\fpbox
\def\abbrevauthors#1{\setbox\aabox=\hbox{\sevenrm\uppercase{#1}}}
\def\abbrevtitle#1{\setbox\atbox=\hbox{\sevenrm\uppercase{#1}}}
\long\def\pag{\beg=\pageno
\def\leftheadline{\noindent\rlap{\nine\folio}\hfil\copy\aabox\hfil}
\def\rightheadline{\noindent\hfill\copy\atbox\hfill\llap{\nine\folio}}
\def\phead{\setbox\fpbox=\hbox{\sevenrm
************************************************}%
\noindent\vbox{\sevenrm\baselineskip9pt\hsize\wd\fpbox%
\centerline{***********************************************}

\centerline{BANACH CENTER PUBLICATIONS, VOLUME **}

\centerline{INSTITUTE OF MATHEMATICS}

\centerline{POLISH ACADEMY OF SCIENCES}

\centerline{WARSZAWA 19**}}\hfill}
\footline{\ifnum\beg=\pageno \hfill\nine[\folio]\hfill\fi}
\headline{\ifnum\beg=\pageno\phead
\else
\ifodd\pageno\rightheadline \else \leftheadline \fi
\fi}}

\newbox\tbox
\newbox\aubox
\newbox\adbox
\newbox\mathbox

\def\title#1{\setbox\tbox=\hbox{\let\\=\cr
\baselineskip14pt\vbox{\twbf\tabskip 0pt plus15cc
\halign to\hsize{\hfil\ignorespaces \uppercase{##}\hfil\cr#1\cr}}}}

\newbox\abbox
\setbox\abbox=\vbox{\vglue18pt}

\def\author#1{\setbox\aubox=\hbox{\let\\=\cr
\nine\baselineskip12pt\vbox{\tabskip 0pt plus15cc
\halign to\hsize{\hfil\ignorespaces \uppercase{\spa{##}}\hfil\cr#1\cr}}}%
\global\setbox\abbox=\vbox{\unvbox\abbox\box\aubox\vskip8pt}}

\def\address#1{\setbox\adbox=\hbox{\let\\=\cr
\nine\baselineskip12pt\vbox{\it\tabskip 0pt plus15cc
\halign to\hsize{\hfil\ignorespaces {##}\hfil\cr#1\cr}}}%
\global\setbox\abbox=\vbox{\unvbox\abbox\box\adbox\vskip16pt}}

\def\mathclass#1{\setbox\mathbox=\hbox{\footnote{}{1991 {\it Mathematics
Subject
Classification}\/: #1}}}

\long\def\maketitlebcp{\pag\unhbox\mathbox
\footnote{}{The paper is in final form and no version
of it will be published elsewhere.}
\vglue7cc
\box\tbox
\box\abbox
\vskip8pt}

\long\def\abstract#1{{\nine{\bf Abstract.}
#1

}}

\def\section#1{\vskip-\lastskip\vskip12pt plus2pt minus2pt
{\bf #1}}

\long\def\th#1#2#3{\vskip-\lastskip\vskip4pt plus2pt
{\sc #1} #2\hskip-\lastskip\ {\it #3}\vskip-\lastskip\vskip4pt plus2pt}

\long\def\defin#1#2{\vskip-\lastskip\vskip4pt plus2pt
{\sc #1} #2 \vskip-\lastskip\vskip4pt plus2pt}

\long\def\remar#1#2{\vskip-\lastskip\vskip4pt plus2pt
\sp{#1} #2\vskip-\lastskip\vskip4pt plus2pt}

\def\Proof{\vskip-\lastskip\vskip4pt plus2pt
\sp{Proo{f.}\ }\ignorespaces}

\def\endproof{\nobreak\kern5pt\nobreak\vrule height4pt width4pt depth0pt
\vskip4pt plus2pt}

\newbox\refbox
\newdimen\refwidth
\long\def\references#1#2{{\nine
\setbox\refbox=\hbox{\nine[#1]}\refwidth\wd\refbox\advance\refwidth by 12pt%
\def\textindent##1{\indent\llap{##1\hskip12pt}\ignorespaces}
\vskip24pt plus4pt minus4pt
\centerline{\bf References}
\vskip12pt plus2pt minus2pt
\parindent=\refwidth
#2

}}

\def\footnoterule{\kern -3pt \hrule width 4cc \kern 2.6pt}

\catcode`@=11
\def\vfootnote#1%
{\insert\footins\bgroup\nine\interlinepenalty\interfootnotelinepenalty%
\splittopskip\ht\strutbox\splitmaxdepth\dp\strutbox\floatingpenalty\@MM%
\leftskip\z@skip\rightskip\z@skip\spaceskip\z@skip\xspaceskip\z@skip%
\textindent{#1}\footstrut\futurelet\next\fo@t}
\catcode`@=12

\font\tt=cmtt10

\mathclass{Primary 58B30; Secondary 81R50 17B37 16W30.}

\abbrevauthors{T. Brzezi\'nski}
\abbrevtitle{Quantum Fibre Bundles}

\title{Quantum Fibre Bundles. An Introduction.}

\author{Tomasz\ Brzezi{\'n}ski}
\address{Institute of Mathematics\\ University of \L\' od\'z\\ ul. Banacha 22\\
90-238 \L\'od\'z, Poland}

\maketitlebcp

\footnote{}{Address after 1st October 1995: University of Cambridge, DAMTP,
Cambridge CB3 9EW, U.K.}

\abstract{An approach to construction of a quantum group
gauge theory based on the quantum group generalisation
of fibre bundles is reviewed.}

\section{1. Introduction and preliminaries.}
\vskip4pt plus2pt

{\bf 1.1.} {\it Introduction}
The algebraic approach to deformation-quantisation involves the
replacing of the
algebras of functions by non-commutative algebras.
In  recent years we have seen a rapid developement of this approach
to quantisation, initiated by  Drinfeld's [15]
realisation of Hopf algebras as deformations of Lie groups. Hopf
algebras are now commonly called quantum groups.
Quantum groups originated in the quantum inverse scattering method developed
by the Petersburg School and applied to quantisation of completely
integrable hamiltonian systems. Nowadays, however, it is believed
that quantisation-deformation and quantum groups in particular  may
be applied to the description of spaces at the Planck scale.
Having this application in mind, it is important to develop
a kind of gauge theory involving quantum groups. Such a theory
was introduced by S. Majid and the author in [6]
in the framework of fibre bundles with quantum structure groups.
In this paper we review the main elements of the
quantum group gauge theory of [6].

The article is organised as follows.
In the remaining part of Section~1
 we give a crash introduction to Hopf algebras
and non-commutative differential geometry. The reader familiar
with these topics may go directly to Section~2, where we describe
elements of the theory of quantum fibre bundles. Then in Section~3
we present gauge theory of such fibre bundles. We conclude
the paper with some remarks on other developments
of quantum group gauge theory and open problems in Section~4.
\vskip4pt plus2pt

{\bf 1.2.} {\it Hopf algebras.}
 A unital algebra $H $ over a field $k$ is called a {\it Hopf algebra} if
there exist  linear maps: a {\it coproduct}  $\Delta :H\to H\otimes H$,  a {\it
counit}  $\epsilon :H\to k$
and an {\it antipode}  $\S: H\to H$  which satisfy the following axioms [26]:

1. $(\Delta\otimes\id)\circ\Delta = (\id\otimes\Delta)\circ\Delta$ ;

2. $(\id\otimes\epsilon)\circ\Delta = (\epsilon\otimes\id)\circ\Delta = \id $;

3. $m\circ(\id\otimes\S) = m\circ(\S\otimes\id) = 1\epsilon$.

Here and in what follows $m$ denotes the mulitpication map.
One should think of a Hopf algebra as a non-commutative generalisation of the
algebra of
regular functions on a group. In this case $\Delta$ corresponds to the group
multiplication and the axiom 1. states the associativity of this
multiplication.
Axiom 2. states the existence of the unit in a group and 3. is the existence
of inverses of group elements, written in a dual form. For this reason Hopf
algebras are also
called {\it quantum groups}.

For a
coproduct we use an explicit
expression $\Delta (a) = a\o\otimes a\t$, where the summation is
implied according to the Sweedler sigma convention [26], i.e.
$a\o\otimes a\t = \sum_{i\in I} a\o^ i\otimes a\t^ i$ for an
index set $I$. We also  use the notation
$$
a\o\tens a\t\tens\cdots\tens a_{(n)} =
(\Delta\tens\underbrace{\id\tens\cdots\tens\id}_{n-2})\circ
\cdots\circ(\Delta\tens\id)\circ\Delta
$$
which describes a multiple action of $\Delta$ on $a\in H$.

A vector space $C$ with a coproduct $\Delta : C\to C\otimes C$
and the counit $\epsilon : C\to k$, satisfying axioms 1. and
2. is called a {\it coalgebra}.

A vector space $V$ is called a {\it right $H$-comodule} if there
exists a linear map $\rho_ R: V\to V\tens H$, called a {\it right
coaction}, such that $(\rho_
R\tens\id)\circ\rho_ R = (\id\tens\Delta)\circ\rho_ R$ and
$(\id\tens \epsilon )\circ\rho_ R = \id$.
We say that a unital algebra $P$ over $k$ is a {\it right H-comodule
algebra} if $P$ is a right $H$-comodule with a coaction
$\Delta_ R :P\to P\otimes H$, and $\Delta_ R$ is an algebra
map.  The algebra structure of $P\tens H$ is that of a tensor product
algebra. For a coaction $\Delta_ R$ we use an explicit notation
$\Delta_ R u = u_{(0)}\otimes u_{(1)}$, where
the summation is also implied. Notice that $u_{(0)}\in
P$ and $u_{(1)}\in H$. If $P$ is a right $H$-comodule so is
$P\otimes P$ with a coaction $\Delta_ R$
$$
\Delta_ R(u\otimes v) = u\z \otimes
v\z\otimes u\o v\o. \eqno{(1)}
$$
If $P$ is a right $H$-comodule algebra then $P^{coH}$ denotes a fixed
point subalgebra of $P$, i.e.
$P^{coH} = \{u\in P :\Delta_ R u = u\otimes 1\}$.
$P^{coH}$ is a subalgebra of $P$ with a natural inclusion $j:
P^{coH} \hookrightarrow P$ which we do not write explicitly later on.

Let $H$ be a Hopf algebra, $B$ be a unital algebra over $k$, and
let $f,g :H\to B$ be linear maps. A {\it convolution product} of
$f$ and $g$ is a linear map $f*g: H\to B$ given by
$(f*g)(a) = f(a\o)g(a\t)$, for any $a\in H$. With respect to
the convolution product, the set of all linear
maps $H\to B$ forms an
associative  algebra with the unit $1\epsilon$. We say that a
linear map $f: H\to B$ is
{\it convolution invertible} if there is a map $f^{-1}: H\to
B$ such that
$f*f^{-1} = f^{-1}*f = 1\epsilon$. The set of all convolution
invertible maps $H\to B$ forms a multiplicative group. Similarly if
$V$ is a right $H$-comodule and $f:V\to B$, $g:H\to B$ are linear maps
then we define a convolution product $f*g:V\to B$ to be $(f*g)(v) =
f(v\z)g(v\o)$.\vskip4pt plus2pt

{\bf 1.3.} {\it Differential structures.}
Let $P$ be a unital algebra over $k$. Denote by $\omup$ the $P$-bimodule
$\ker m$,
where $m: P\otimes P \to P$ is a multiplication map. Let $\du : P\to \omup$ be
a
linear map
$$
\du u = 1\otimes u - u\otimes 1. \eqno{(2)}
$$
 It can be easily checked that $\du$ is a differential, known as the Karoubi
differential. We call the pair $(\omup, \du)$ the {\it universal differential
structure on P} [19, 20]. $\omup$ should be understood as a bimodule
of 1-forms. We say that $ (\Omega^ 1(P) ,\d)$
is a {\it first order differential calculus} on $P$ if there exists
a subbimodule $\CN\subset\omup$ such that $ \omp = \omup /\CN$ and
$\d = \pi\circ\du$, where $\pi: \omup\to\omp$ is a canonical projection.
It is then said that $(\omp, \d)$  is generated by $\CN$. Let a differential
structure $(\Omega^ 1(H), \d)$ on a Hopf algebra $H$ be generated by
$\CN \subset\Omega^ 1H$. We say that  $(\Omega^ 1(H), \d)$
is a {\it bicovariant differential calculus} [29] if there exists a unique
right ideal $\CQ\subset\ker\epsilon$ such that $H\otimes\CQ =\kappa(\CN)$,
where $\kappa :H\otimes H\to H\otimes H$, $\kappa : a\otimes b \mapsto
ab\o\otimes b\t$, and $\ad(\CQ) \subset \CQ\otimes H$, where $\ad : H\to
H\otimes H$
is a right adjoint coaction
$$
\ad :a\mapsto a\t\otimes (\S a\o)a_{(3)}.\eqno{(3)}
$$
The universal
differential envelope is the unique differential algebra $(\Omega P, \d)$
containing $(\omup, \du)$ as its 1-st order part.

\section{2. Fibre bundles.}
 In this section we report the basic elements of the theory of
quantum fibre bundles of S. Majid and the author
[6]. The detailed analysis of quantum
group gauge theory on classical spaces may be found in
[7]. All the algebras are over a field $k$ of complex
or real numbers. Except for Section~2.4 and
Example~3.1.4
 we work with the universal differential structure.
\vskip4pt plus2pt

{\bf 2.1.} {\it Quantum principal bundles.}
Let $H$ be a Hopf algebra, $P$ a right $H$-comodule algebra
with a coaction
$\Delta_ R :P\to P\otimes H$. We define a canonical map $\chi :P\otimes P \to
P\otimes H$,
$$
\chi = (m\otimes \id)\circ (\id \otimes \Delta_ R).
\eqno{(4)}
$$
Explicitly, $\chi (u\otimes v) = uv\z\otimes
v\o, $ for any $u,v\in P$. We say that the coaction
$\Delta_ R$ is
{\it free} if $\chi$ is a surjection and it is {\it exact}
if $\ker \chi = P(\d P^{coH})P,$ where $\d$ denotes the
universal differential (2) and $P^{coH}$ is a fixed point subalgebra
of $P$.  We denote $P(\d P^{coH})P$ by $\omup_{\rm hor}$ and call
its elements {\it horizontal forms}. Although the freeness and exactness
conditions are algebraic in this formulation one should notice
that in fact the latter one is a condition on differential
structures on $P$ and $P^{coH}$. This becomes clear
in Section~2.4.
The map $\chi\mid_\omup$ has a natural geometric
interpretation as a dual to the map $ {\cal G}\to T_ uX$, which to
each element of the Lie algebra $\cal G$ of a group $G$
associates a fundamental vector field on a manifold  $X$ on which $G$
acts.
\defin{Definition}{2.1.1 Let $H$ be a Hopf algebra, $(P ,\Delta_{R})$ be a
right $H$-comodule
algebra and let  $B = P^{coH}$. We say that $P(B,H)$ is a {\it quantum
principal bundle} within the differential
envelope, with a structure quantum group $H$ and a base $B$ if the
coaction $\Delta_ R$ is free and exact.}

This definition reproduces the classical situation (but in a dual
language)  in which a group $G$ acts freely on a total space $X$ from
right, and a base  manifold $M$ is defined as $M = X/G$. The freeness
of the action of $G$ on  $X$ means that a map
$X\times G \to X\times X$, $(u,g)\mapsto (u,ug)
$ is an inclusion. In the classical situation and the  commutative
differential structure the exactness follows from the
freeness. This
is no longer true in a non-commutative extension.

The notion of a quantum principal bundle is strictly related to
the theory of  algebraic extensions [25] since $P(B,H)$
is a Hopf-Galois extension
of $B$ to $P$ by a Hopf algebra $H$.
Yet another way of defining of a quantum principal bundle makes use
of the notion of a {\it translation map}, which proves very useful
in analysis of the structure of quantum bundles [4].
\th{Proposition}{2.1.2.}{
Let $H$ be a Hopf algebra, $P$ a right
$H$-comodule algebra and $B = P^{co H}$.
Assume that the coaction $\Delta_ R$ is free.
Then $P(B,H)$ is a quantum principal bundle iff there exists a linear  map
$\tau:H\to
P\otimes{}_ B P$,  given by $\tau(a) = \sum_{i\in I}u_
i\otimes{}_ B v_ i$, where $\sum_{i\in I}u_ i\otimes v_ i \in
\chi^{-1}(1\otimes a)$. The map $\tau$  is called a {\it translation map}.}

A translation map is a well-known object in the classical bundle theory
[18, Definition~2.1].
Classically, if $X$ is a manifold on which a Lie group $G$ acts
freely then the translation map $\hat{\tau}:X\times{}_ MX\to G$, where
$M=X/G$, is defined by $u\hat{\tau}(u,v)=v$.
Dualising this construction we arrive immediately at
the map $\tau$ above.
\vskip4pt plus2pt

{\bf 2.2.} {\it Examples of quantum principal bundles}

\remar{Example\ {2.2.1.}\ }
{{\it A trivial quantum principal bundle.}
 Let $H$ be a Hopf algebra, $P$  a right $H$-comodule
algebra  and $B = P^{coH}$.
Assume there is a convolution invertible map $\Phi :
H\rightarrow P$ such that
$
\Delta_ R\Phi = (\Phi\otimes\id)\Delta ,\quad \Phi(1)=1 ,
$
i.e.~$\Phi$ is an intertwiner. Then $P(B,H)$ is a quantum principal
bundle  called a {\it trivial quantum principal bundle} and denoted by
$P(B,H,\Phi)$.  The word  {\it trivial} refers to the fact that $P
\cong B\otimes H$ as  vector spaces with an isomorphism $\Theta_\Phi
: P\to B\otimes H$, $\Theta_\Phi:
u\mapsto u\z\Phi^{-1}(u\o)\otimes u\t$. Moreover, as algebras
$P\cong B_\Phi\# H$, where
${}_\Phi\#$ denotes a crossed product [1], with the
isomorphism $\Theta_\Phi$
above. Explicitly, the product in
$B_\Phi\# H$ is given by
$$
(b_ 1\otimes a^ 1)(b_ 2\otimes a ^ 2) =  b^ 1\Phi( a^1\o) b_2
\Phi(a^2\o)\Phi^{-1}(a^1\t a^2\t)\otimes a^1_{(3)} a^2_{(3)} .
$$
Such an algebra $P$ is also known as a {\it cleft extension} of $B$
[27, 14].

The  map $\tau = (\Phi^{-1}\otimes_ B\Phi)\circ\Delta$ is a translation
map in $P(B,H,\Phi)$.}

For a trivial quantum principal bundle $P(B,H,\Phi)$ we define a
{\it gauge transformation}  as a convolution invertible  map $\gamma:
H\to B$ such that $\gamma(1) =1$. The set of all gauge transformations
of $P(B,H,\Phi)$ forms a group with respect to the
convolution product. This group is denoted by $\CH(B)$. Gauge
transformations relate different trivialisations of $P(B,H,\Phi)$:
$\Psi :H\to P$ is a trivialisation of $P(B,H,\Phi)$ iff there exists $\gamma
\in \CH(B)$ such that $\Psi = \gamma*\Phi$. They also
have a clear meaning in the theory of crossed products. The
following proposition is a special case of the result of Doi [13]
(see also [21, Proposition 4.2]).
\th{Proposition}{2.2.2.}{
Let $P(B,H,\Phi)$ be a trivial quantum principal bundle. Let for any
trvialisation $\Psi$ of $P(B,H,\Phi)$, $\Theta_\Psi :
B_\Psi\rtimes H\to B_\Phi\rtimes H$ be a crossed product algebra
isomorphism such that $\Theta_\Psi\mid_ B = \id$ and $\Delta_ R\Theta
_\Psi = (\Theta_\Psi\otimes\id)\Delta_ R$. Then there is a
bijective correspondence between all isomorphisms $\Theta_\Psi$
corresponding to all trivialisations $\Psi$ and the gauge transformations
of $P(B,H,\Phi)$.}

\remar{Example\ {2.2.3.}\ }{
{\it Quantum principal bundle on a quantum homogeneous space.}
 Let $H$ and $P$ be Hopf algebras. Assume, there is a Hopf algebra
projection  $\pi : P\rightarrow H$. Define a right coaction of $H$ on
$P$ by
$
\Delta_ R = (\id\otimes\pi)\Delta : P\rightarrow P\otimes H.
$
Then $B = P^{coH}$ is a {\it quantum quotient space}, a special
case of a {\it quantum homogeneous space}.
Assume that $\ker\pi \subset m\circ (\ker\pi\mid_ B\otimes P)$. Then
$P(B,H)$ is a quantum principal bundle within the differential
envelope. This bundle is denoted by $P(B,H,\pi)$.

The translation map $\tau: H\to
P\otimes{}_ B P$ in $P(B,H,\pi)$ is given by
$\tau(a) = \S u\o\otimes{}_ B u\t$,
where $u\in\pi^{-1}(a)$.}

A large number of examples of quantum bundles on quantum homogeneous
spaces has been found in [22]. The simplest
and probably the most fundamental one is
\remar{Example\ {2.2.4.}\ }
{{\it The quantum Hopf fibration \rm [6, Section~5.2]}. The
total space of this bundle is the quantum group $SU_
q(2)$, as an algebra generated by the identity and a matrix
$T = (t_{ij}) = \pmatrix{\alpha & \beta \cr \gamma &\delta}$, subject
to the homogeneous relations
$$
\alpha\beta = q\beta \alpha , \quad \alpha\gamma = q
\gamma\alpha , \quad \alpha \delta = \delta\alpha +
(q-q^{-1})\beta\gamma ,\quad
\beta\gamma = \gamma\beta, \quad \beta\delta = q\delta \beta
,\quad \gamma \delta = q \delta \gamma ,
$$
and a determinant relation
$\alpha\delta-q\beta\gamma=1$, $q\in k^*$. $SU_ q(2)$ has a matrix quantum
group
structure,
$$
\Delta t_{ij} = \sum_{k =1}^ 2 t_{ik}\tens t_{kj}, \quad
\eps(t_{ij})=\delta_{ij},\quad \S T = \pmatrix{\delta &-q^{-1}
\beta \cr -q\gamma & \alpha}.
$$
The structure quantum group of the quantum Hopf bundle is an algebra
of functions on $U(1)$, i.e. the algebra $k[Z,Z^{-1}]$ of formal
power series in $Z$ and $Z^{-1}$, where $Z^{-1}$ is an inverse of
$Z$. It has a standard Hopf algebra structure
$$
\Delta Z^{\pm 1} = Z^{\pm 1}\tens Z^{\pm 1} , \quad
\eps(Z^{\pm 1}) = 1, \quad
\S Z^{\pm 1} = Z^{\mp 1} .
$$
There is a Hopf algebra projection $\pi: SU_ q(2) \to
k[Z,Z^{-1}]$,
$$
\pi : \pmatrix{\alpha & \beta \cr \gamma &\delta} \mapsto \pmatrix{Z &
0\cr 0 & Z^{-1}},
$$
which defines a right coaction $\Delta_ R : SU_ q(2) \to SU_
q(2) \tens k[Z,Z^{-1}]$ by $\Delta_ R =(\id\tens\pi)\circ\Delta$.
Finally $S_ q^ 2\subset SU_ q(2)$ is a quantum two-sphere
[24], defined as a fixed point subalgebra, $S_ q^ 2 = SU_
q(2)^{cok[Z,Z^{-1}]}$. $S_ q^ 2$ is generated by $\{1, b_ - =
\alpha\beta , b_ + =\gamma\delta ,b_ 3 = \alpha\delta\}$ and the
algebraic relations in $S_ q^ 2$ may be deduced from those in
$SU_ q(2)$.

It was shown in [6] that $SU_
q(2)(S_ q^ 2,k[Z,Z^{-1}],\pi)$ is a non-trivial quantum
principal bundle over the homogeneous space.}

The other examples of quantum principal bundles constructed in [22]
include:
$$
U_ q(n)(S_ q^{2n-1}, U_ q(n-1),\pi),$$
$$
SU_ q(n)(S_ q
^{2n-1}, SU_ q(n-1),\pi),$$
$$SU_ q(n)(\C\P_ q^{n-1}, U_ q(n-1),
\pi),$$
$$
U_ q(n)(G_ k(\C^ n_ q), U_ q(k)\otimes U_ q(n-k), \pi),
$$
where $G_ k(\C^ n_ q)$ is a quantum Grassmannian.

\remar{Remark\ {2.2.5.}\ }{The quantum sphere $S^2_q$ considered
in Example~2.2.4. is the special case of the most general
quantum sphere $S^2_q(\mu ,\nu)$, where $\mu\neq\nu$ are
real parameters such that $\mu \nu \geq 0$ (see [24] for
details). Precisely $S^2_q = S^2_q(1,0)$. It can
be shown that $S^2_q$ is the only quantum sphere
which can be interpreted as a quotient space of $SU_q(2)$
by $k[Z,Z^{-1}]$ in the sense
of Example~2.2.3. It turns out, however, that $S^2_q(\mu ,\nu)$
may be veiwed as a quotient space of $SU_q(2)$ by
a {\it coalgebra} $C = SU_q(2)/J$, where $J$ is a right
ideal in $SU_q(2)$ generated by
$$
p(q\alpha^2 - \beta^2) + \alpha\beta -pq, \quad
p(q\gamma^2 - \delta^2) + \gamma\delta + p, \quad
p(q\alpha\gamma - \beta\delta) + q\beta\gamma ,
$$
where $p = \sqrt{\mu\nu}/(\mu-\nu)$ [5].  Precisely
$$
S^2_q(\mu,\nu) = \{u\in SU_q(2);\;\; u\o\otimes\pi(u\t) =
u\otimes\pi( 1)\},
$$
where $\pi:SU_q(2)\to C$ is the canonical surjection.
It can be shown that the vector space $C$ is spanned
by $1=\pi(1)$, $x_ n = \pi(\alpha^n)$ and
$y_n = \pi(\delta^n)$ (cf. definition of $\pi$ in Example~2.2.4).

One would like to view $SU_q(2)$ as a total space of a quantum
principal bundle over $S_q^2(\mu, \nu)$similarly as in Example~2.2.4. Since
$C$ is not a Hopf algebra one needs to generalise the notion
of a bundle. In [5] we proposed the following generalisation
of Definition~2.1.1. Let $C$ be a coalgebra and let $P$ be an algebra
and a right $C$-comodule. Assume that there
is an action $\rho : P\otimes C\otimes P \to P\otimes
C$ of $P$ on $P\otimes C$ and an element $ 1\in C$
such that $\Delta_R\circ m = \rho\circ(\Delta_R\otimes\id)$ and
for any $u, v\in P$, $\rho(u\otimes  1 , v)
= \chi(u\otimes v)$. Then $B = \{ u\in P;\;\; \Delta_R u =u\otimes 1\}$
is a subalgebra of $P$, and we say that $P(B,C,\rho)$
is a {\it quantum $\rho$-principal bundle} over $B$
if the coaction $\Delta_R$ is free and exact.

In the above example of the quantum sphere $S^2_q(\mu,\nu)$ the action
$\rho$ is given by $\rho(u\otimes c, v) = uv\o\otimes
\rho_0(c,v\t)$, where $\rho_0$ is a natural right
action of $SU_q(2)$ on $C$.}
\vskip4pt plus2pt

{\bf 2.3.} {\it Quantum associated bundles}
\defin{Definition} {2.3.1.
\rm Let $P(B,H)$ be a quantum principal bundle and let $V$ be a right
$H^{\rm op}$-comodule algebra, where $H^{\rm op}$ denotes the algebra
which is isomorphic to $H$ as a vector space but has an opposite
product, with coaction $\rho_{R} : V \rightarrow V
\otimes H$. The space $P \otimes V$ is naturally endowed with a right
$H$-comodule structure  $\Delta_{E} : P \otimes V
\rightarrow P \otimes V \otimes H$ given by
$\Delta_{E} (u \otimes v) =  u\z \otimes
v\z \otimes u\o v\o$
for any $u \in P$ and $v \in V$. We say that the fixed point subalgebra
$E$ of $P\tens H$ with respect to $\Delta_ E$
is a  {\it quantum fibre bundle associated to $P(B,H)$} over $B$ with
structure quantum group $H$ and standard fibre $V$. We denote it by $E =
E(B,V,H)$.}

It can be easily shown that $B$ is a subalgebra of $E$ with the inclusion
$j_ E = b\tens 1$. The inclusion $j_ E$ provides $E$ with the
structure of a left $B$-module.
\remar{Example\ {2.3.2.}\ }
{ Let $P(B,H,\Phi)$ be a trivial quantum principal bundle and let $V$ be
as in Definition~2.3.1. Assume also that $H$
has a bijective antipode. The associated bundle $E(B,V,H)$ is called
a {\it trivial quantum fibre bundle}. Trivialisation $\Phi :H\to P$ induces a
map $\Phi_ E :V\to E$,
$
\Phi_{E} (v) = \sum \Phi (\S^{-1} v\o) \otimes
v\z$
which allows one to identify $E$ with $B\tens V$ as vector spaces via
the linear isomorphism $b\tens v\mapsto b\Phi_ E(v)$. As an algebra,
$E$ is isomorphic to a certain crossed product algebra $B\rtimes V$
[2]. }

The following proposition shows that a
quantum principal bundle is a fibre bundle associated to itself.
\th{Proposition}{2.3.3.} {
A quantum principal bundle $P(B,H)$ is a fibre bundle associated to
$P(B,H)$ with the fibre which is isomorphic to $H$ as an algebra and
with the coaction $\rho_ R=(\id\otimes \S)\circ\Delta'$, where
$\Delta'$ denotes the opposite coproduct, $\Delta'(a) = a\t\tens a\o$,
for any $a\in H$.}

{}From the point of view of a gauge theory it is important to consider
cross-sections
of a vector bundle. In this algebraic setting a cross-section is defined as
follows
\defin{Definition}{ 2.3.4.
 Let $E(B,V,H)$ be a quantum fibre bundle associated to a quantum
principal bundle $P(B,H)$. A left $B$-module map
$s: E\to B$ such that $s(1) = 1$ is called a {\it cross section} of
$E(B,V,H)$. The set of cross sections of $E(B,V,H)$ is denoted by
$\Gamma(E)$. }
\th{Lemma}{ 2.3.5.} {
If $s:E\to B$ is a cross section of a quantum fibre bundle $E(B,V,A)$  then
$s\circ j_ E = \id$.}

The result of trivial Lemma~2.3.5. justifies the term
cross section used in Definition~2.3.4. We  remark that the
definition of a cross section of a quantum fibre bundle analogous to the
one we use
here was first proposed  in [17].
We analyse cross-sections more closely in Section~3.3.
\vskip4pt plus2pt

{\bf 2.4.} {\it Quantum principal bundles with general differential
structures.}
The detailed analysis of quantum principal bundles with general
differential structures goes far beyond the scope of this
paper. Here we give only a definition of a quantum principal bundle
with general differential structure. We refer the interested reader to the
fundamental
paper [6]. More explicit exposition may be
also found in [2].

Let $(\omp ,\d)$ be a first order differential calculus on a right $H$-comodule
algebra $P$ generated
by $\CN\subset\omup$ and let $(\Omega^ 1(H), \d)$ be a bicovariant
differential structure on $H$ generated by the right ideal $\CQ\subset
\ker\epsilon$. We say that differential structures $(\omp, \d)$
and $(\Omega^ 1(H),\d)$ {\it agree} with each other if
$\Delta_ R(\CN)\subset\CN\otimes H$ , where $\Delta_ R$ is given
by (1),  and $\chi(\CN)\subset P\otimes\CQ$.
If differential structures on $P$ and $H$ agree we can define
a map $\chi_\CN :\omp\to P\otimes\ker\epsilon /\CQ$ as follows.
Let $\pi_\CN :\omup\to\omp$ and $\pi_\CQ :\ker\epsilon\to
\ker\epsilon /\CQ$ be canonical projections. Then for any $\rho\in\omp$
take any $\rho_ U\in\pi_\CN^{-1}(\rho)$ and define
$\chi_\CN(\rho) = (\id\otimes\pi_\CQ)\circ\chi(\rho_ U)$, where
$\chi$ is a canonical map (4). We say that the coaction
$\Delta_ R:P\to P\otimes H$ is {\it exact} with respect to differential
structures generated by $\CN$ and $\CQ$ if $\ker\chi_\CN =
P\Omega^ 1(P^{coH}) P$. Finally we define a quantum principal
bundle with $P(B,H)$ with differential structure generated by $\CN$
and $\CQ$ if the coaction $\Delta_ R$ is free and exact with respect
to this structure.

\section{3. Gauge Theory.}
In this section we analyse more closely the
structure of quantum bundles. We introduce
the formalism of connections and take a closer look at cross sections and
gauge transformations in general (non-trivial) quantum bundles.
\vskip4pt plus2pt

{\bf 3.1. }{\it Connections = gauge fields.}
{}From the point of view of gauge theories principal connections
are the gauge fields. In the definition of a principal connection an important
r\^ole is played by a right adjoint coaction of $H$ on itself (3).
Since  $\ad (\ker \epsilon)\subset\ker\epsilon\otimes H$,
 we can define a coaction $\Delta_ R
:P\otimes\ker\epsilon \to P\otimes\ker\epsilon\otimes H$ by
$
\Delta_ R (u\otimes a) = u\z\otimes a\t\otimes
u\o(\S a\o)a_{(3)}.
$The canonical map $\chi :\omup\to P\otimes\ker\epsilon$
is equivariant, i.e.
$
\Delta_ R \chi = (\chi\otimes \id )\Delta_
R,
$
where $\Delta_ R$ on $\omup$ is given by (1).
{}From the definition of a quantum principal bundle we deduce
that the following sequence
$$
0\to\Gamma_{hor}\buildrel{j}\over{\to}\omup\buildrel{\chi}\over\to
P\otimes\ker\epsilon\to 0
$$
is an exact sequence of equivariant maps. A connection in
$P(B,H)$ is a right-invariant splitting of this sequence. In
other words, if there is a map
$
\sigma
:P\otimes\ker\epsilon\to\omup
$
such that $\Delta_
R\sigma = (\sigma\otimes\id)\Delta_{R}$ and
$\chi\circ\sigma = \id$, then a connection in
$P(B,H)$ is identified with a linear projection $\Pi: \omup\to\omup$, $\Pi =
\sigma\circ\chi\mid_\omup$. Obviously, $\Delta_ R\Pi =
(\Pi\otimes\id)\Delta_ R$. The connection $\Pi$ is {\it strong} if and only if
$(\id -\Pi)\d P \subset \Omega^ 1BP$, [17].

We denote $\omup_{\rm ver} = {\rm Im}\Pi$. Every
$\alpha\in\omup_{\rm ver}$ is said to be a {\it vertical
1-form}. If there is a connection in $P(B,H)$, then
$\omup = \omup_{\rm hor}\oplus\omup_{\rm ver}$.

Next we define a map $\omega :H\to \omup$, by
$$
\omega(a) = \sigma(1\otimes (a-\epsilon(a))).
$$
The map $\omega$ is called a {\it connection 1-form} of the
connection $\Pi$.

\th{Theorem} {3.1.1.} {Let $P(B,H)$ be a quantum principal bundle and let
$\Pi$
be a connection in
$P(B,H)$. A connection form $\omega$ has the following properties:

1. $\omega(1) = 0$;

2. $\forall\; a\in H, \quad\chi\omega (a) = 1\otimes
(a- \epsilon (a))$;

3. $\Delta_{R} \circ\omega = (\omega\otimes\id)\circ\ad$.

Conversely, if $\omega :H\to \omup$ is a linear map
obeying 1-3, then $\Pi = m\circ (\id\otimes\omega)\chi\mid_\omup$ is
a connection with a connection 1-form $\omega$.}

Having a connection $\Pi$ in a quantum principal bundle  $P(B,H)$ one can
define
the horizontal projection as a complimentary part of $\Pi$, and a covariant
derivative
as a horizontal part of $\d$ (for details see [6]). As a result one
defines a curvature of a strong connection $\omega$ as $F =
\d\omega+\omega*\omega$
[17].
\remar{Example\ {3.1.2.}\ }{
{\it Strong connection in a trvial bundle.}
 Let $P(B,H,\Phi)$ be a trivial quantum principal bundle as
before, and let $\beta :H\to \Omega^ 1B$ be any linear
map such that $\beta(1)=0$. Then the map
$
\omega = \Phi^{-1}*\beta*\Phi +\Phi^{-1}*\d\Phi
$
is a connection 1-form in $P(B,H,\Phi)$. Its curvature is easily computed
to be $F = \Phi^{-1}*(\d\beta +\beta*\beta)*\Phi$.}
\remar{Example\ {3.1.3.}\ }
{{\it Canonical connection. } Let $P(B,H,\pi)$ be a quantum principal bundle
over the
homogeneous space $B$ as described in
Example~2.2.3.  Assume, there is an algebra
inclusion $i :H\hookrightarrow P$ such that $\pi\circ i =\id$,
$\epsilon_ P(i(a)) = \epsilon_ H (a)$, for any $a\in H$ and such
that
$
(\id\otimes\pi) \ad i = (i\otimes \id) \ad .
$
Then the map $\omega(a) = \S i(a)\o\d i(a)\t$
is a connection 1-form in $P(B,H,\pi)$. This connection
is strong if $i$ is an
intertwiner for the right coaction [2, Lemma~5.5.5].}

\remar{Example\ {3.1.4.}\ }
{{\it The Dirac $q$-monopole.}
 Consider the quantum Hopf fibration of Example~2.2.4.
Let a differential structure $(\Omega^ 1(SU_ q(2)), \d)$ be given by the
3D caluclus of Woronowicz [28]. $\Omega^ 1(SU_ q(2))$
is generated by  the
forms $\omega^ 0 = \delta\d\beta - q^{-1}\beta\d\delta$,
$\omega^ 1 = \delta\d\alpha - q^{-1}\beta\d\gamma$,
$\omega^ 2 = \gamma\d\alpha - q^{-1}\alpha\d\gamma$
and the
relations
$$
\omega^0\alpha  =  q^{-1}\alpha\omega^0 , \quad \omega^0\beta =
q\beta\omega^0 ,\quad
\omega^1\alpha  =  q^{-2}\alpha\omega^1 ,
$$
$$
\omega^1\beta  =
q^2\beta\omega^1, \quad
\omega^2\alpha  =  q^{-1}\alpha\omega^2 , \quad \omega^2\beta =
q\beta\omega^2 .
$$
The remaining relations can be obtained by the replacement
$\alpha\rightarrow\gamma$, $\beta\rightarrow\delta$. One can show
that $SU_ q(2)(S_ q^ 2, k[Z,Z^{-1}], \pi)$ is a quantum
principal bundle with this differential structure. We define
the connection one form $\omega:k[Z,Z^{-1}]\to\Omega^{1}(SU_ q(2))$
by
$$
\omega(Z^ n) = {{q^{-2n} - 1}\over{q^{-2} -1}}\omega^ 1.
$$
In [6] it has been shown that $\omega$ is a canonical connection in
$SU_ q(2)(S_ q^ 2, k[Z,Z^{-1}], \pi) $ which reduces to
the Dirac monopole of charge 1 [16] when $q\to 1$. The curvature
of $\omega$ is $F(Z^ n) = {{q^{-2n} - 1}\over{q^{-2} -1}}\omega^ 0
\wedge\omega^ 2$. The q-deformed Dirac monopole of any charge is
 discussed
in [12].}
\vskip4pt plus2pt

{\bf 3.2.} {\it Cross sections = matter fields}
In this section we use the notion of a translation map in a quantum
principal bundle $P(B,A)$ to identify cross sections of a quantum
fibre bundle $E(B,V,A)$  with equivariant maps $V\to P$. In gauge theories
such maps play a role of matter fields. Recall that a
linear map $\phi:V\to P$ is said to be equivariant if $\Delta_ R
\phi = (\phi\tens\id)\rho_ R$, where $\rho_ R$ is a right coaction
of $A$ on $V$. In particular, our identification
implies that a quantum principal bundle is trivial if it admits a
cross section which is an algebra map.

\th{Theorem} {3.2.1} {
Let $H$ be a Hopf algebra with a bijective antipode. Cross sections of
a quantum fibre bundle $E(B,V,H)$ associated to a
quantum principal bundle $P(B,H)$ are in bijective
correspondence
with  equivariant maps $\phi :V\to P$ such that $\phi(1) = 1$.}
\Proof A map $\phi :V \to P$ induces a
cross section $s$ of $E(B,V,H)$,  by $s =
m\circ(\id\tens\phi)$. Conversely, for any $s\in\Gamma(E)$ we define a
map $\phi: V\to P$ by
$$
\phi : v\mapsto \tauo(\S^{-1}v\o)s(\taut(\S^{-1}v\o)\tens v\z),
\eqno{(5)}
$$
where $\tau(a) = \tauo(a)\tens{}_ B\taut(a)$ is a translation map in
$P(B,H)$, and then use properties of a translation map to prove that $\phi$ has
the
required properties and that the correspondence $\theta: \phi\mapsto s$ is
bijective.
\endproof

\remar{Example\ {3.2.2.}\ }
{ Let $E(B,V,H)$ be a quantum fibre bundle associated to a trivial
quantum principal bundle $P(B,H,\Phi)$ as described in
Example~2.3.2. In this case every element of
$E$ has the from $\sum_{i\in I}b_ i\Phi_ E(v_ i)$ for some
$b_ i\in B$ and $v_ i \in V$, and the bijection $\theta$ of the
proof of Theorem~3.2.1 reads
$$
\theta(\phi)(\sum_{i\in I}b_ i\Phi_ E(v_ i)) = \sum_{i\in
I}b_ i\Phi(\S^{-1}v_ i\o)\phi(v_ i\z),
$$
for any equivariant $\phi:V\to P$. The inverse of $\theta$ associates
an equivariant map $\theta^{-1}(s):V\to P$,
$$
\theta^{-1}(s)(v) = \Phi^{-1}(\S^{-1}v\o)s(\Phi_ E(v\z))
$$
to any $s\in\Gamma(E)$. Notice that the map $\theta^{-1}(s)$
obtained in this way is different from the equivariant map $\phi$
discussed in [6, Proposition A6].}

\th{Corollary }{3.2.3} {
Cross sections $s:P\to B$ of a quantum principal bundle $P(B,H)$ are in
bijective correspondence
with the maps $\phi :H\to P$ such that $\Delta_ R\phi =
(\phi\tens\S)\Delta'$ and $\phi(1) = 1$.}

Note that in Corollary~3.2.3. we do not
need the invertibility of $\S$, but
 if $H$   has a  bijective antipode $\S$, the
sections of
a quantum principal bundle $P(B,H)$ are in one-to-one correspondence
with the maps $\psi:H\to P$ such that $\psi(1) = 1$ and $\Delta_
R\circ\psi = (\id\tens\psi)\circ\Delta$. We simply need to define
$\psi=\phi\circ\S^{-1}$, where $\phi$ is given by
Corollary~3.2.3.
\th{Proposition }{3.2.4 }{
Any trivial quantum principal bundle $P(B,H,\Phi)$ admits a section.
Conversely, if a bundle $P(B,H)$ admits a section which is an algebra
map then $P(B,H)$ is trivial with the total space $P$ isomorphic to
$B\tens H$ as an algebra.}

\Proof A convolution inverse of a trivialisation $\Phi$ of a trivial quantum
principal bundle
$P(B,H,\Phi)$ satisfies the assumptions of
Corollary~3.2.3, hence $s=\id*\Phi^{-1}$ is a
section of $P(B,H,\Phi)$. Conversely, assume that an algebra map
$s:P\to B$ is a section of $P(B,H)$. Clearly, $s$ is a $B$-bimodule
map, hence we can define a linear map $\Phi:H\to P$, $\Phi
=m\circ(s\otimes{}_ B\id)\circ\tau$. One then shows that $\Phi$ is a
trivialisation and
$\tilde{\theta}(s)$ constructed in
Corollary~3.2.3 is its convolution inverse. \endproof
\remar{Remark\ {3.2.5.}\ }
{ We would like to emphasise that the existence of a cross section of a
quantum principal bundle does not necessarily imply that the bundle
is trivial. As an example of a non-trivial quantum principal bundle
admitting a cross section we consider  the quantum Hopf fibration of
Example~2.2.4. We consider a linear map
$\phi : k[Z,Z^{-1}]\to SU_ q(2)$, given by
$$
\phi(1) = 1,\qquad \phi(Z^ n) = \delta^ n, \qquad \phi(Z^{-n})
=\alpha^ n,
$$
for any positive integer $n$. The map $\phi$ satisfies the hypothesis
of Corollary~3.2.3, hence it induces a cross section $s:
SU_ q(2) \to S_ q^ 2$, $s: u\mapsto u\o\phi(\pi(u\t))$ but
$s$ is not an algebra map since, for example,
$
s(\alpha\beta) = b_ -\neq q^{-1}b_ 3 b_ - = s(\alpha)s(\beta).
$}
\vskip4pt plus2pt
{\bf 3.3.} {\it Vertical automorphisms = gauge transformations}
\defin{Definition}{ 3.3.1.
 Let $P(B,H)$ be a quantum principal bundle. Any left $B$-module
automorphism $\CF : P\to P$ such that $\CF (1) =1$ and $\Delta_ R
\CF = (\CF\otimes \id)\Delta_ R$ is called a {\it vertical
automorphism} of the
bundle $P(B,H)$. The set of all vertical automorphisms of $P(B,H)$ is
denoted by $Aut_ B(P)$.}

Elements of $Aut_ B(P)$ preserve both the base space $B$ and the
action of the structure quantum group $H$ of a quantum principal
bundle $P(B,H)$. $ Aut_ B(P)$ can be equipped with a multiplicative
group structure $\cdot :(\CF_ 1, \CF_ 2)\mapsto \CF_
2\circ\CF_ 1$. Vertical automorphisms are often called gauge
transformations and $ Aut_ B(P)$ is termed a gauge group.
\th{Proposition} {3.3.2.} {
 Vertical automorphisms
of a quantum principal bundle $P(B,H)$ are in bijective correspondence
with convolution invertible
maps $f:H\to P$ such that $f(1)=1$ and $\Delta_ Rf = (f\otimes
\id)\ad$.}

\Proof If $f$ is a map satisfying the hypothesis of the proposition. then
  $\CF =\id*f$. Conversely, for any $\CF\in Aut_ B(P)$ a map
$f: A\to P$,
$f =m\circ(\id\tens{}_ B\CF)\circ\tau$ ,
where $\tau$ is a translation map has all the required properties. \endproof

Maps $f:H\to P$ form a
group with respect to the convolution product. This group is
denoted by $\CH(P)$. There is an action of $\CH(P)$ on
the space of connection one-forms in $P(B,A)$ given by
$(\omega, f)\mapsto \omega^ f = f^{-1}*\omega*f + f^{-1}*\d f$.
The connection one-form $\omega^ f$ is called a gauge transformation
of $\omega$. If $\omega$ is strong so is its gauge transformation.
Gauge transformation of such $\omega$ induces the gauge
transformation of its curvature $F\mapsto f^{-1}*F*f$.
Similarly there is an action of $\CH(P)$ on $\Gamma(E)$ viewed
as equivariant maps $\phi: V\to P$ by Theorem~3.2.1.,
given by $(\phi, f)\mapsto\phi^ f =\phi*f$. These are the transformation
properties of the fields in quantum group gauge theories.

Proposition~3.3.2. implies the following:
\th{Corollary }{ 3.3.3. } {
For a quantum principal bundle $P(B,H)$, $Aut_ B(P)\cong\CH(P)$ as
multiplicative groups.}
\th{Theorem }{ 3.3.4. }{
Let $P(B,H,\Phi)$ be a trivial quantum principal bundle. Then the
groups $Aut_ B(P)$, $\CH(P)$, and the gauge group
$\CH(B)$ are isomorphic to each other.}

Therefore
Theorem~3.3.4. allows one to interpret a vertical
automorphism of a (locally) trivial quantum principal bundle as
a change of local variables and truly as a gauge transformation of a
trivial quantum principal bundle.

\section{4. Conclusions and open problems}
In this paper we reviewed basic properties
of quantum fibre bundles introduced in [6]. There is
a number of constructions, already present in the literature,
 that we have not described in here.
 For example, locally trivial quantum principal bundles,
defined in [6] were developed
by M. Pflaum in [23], using the methods of the sheaf
theory. A very interesting example of the Yang-Mills theory
in quantum bundles was constructed by P. Hajac in [17].
The example considered in [17]
 belongs to the interface of the theory
described here and the Connes-Rieffel Yang-Mills
theory [11], and points to the very important problem of finding
the relationship between the quantum group gauge theory
and Connes' non-commutative geometry [9].

There is also a number of challenging problems
that need to be solved in order to obtain a full understanding
of quantum group gauge theories. For example, in this article
we restricted our discussion only to gauge transformations
of bundles with the universal differential structure. The theory
of gauge transformations of bundles with general differential
structures is not yet known. In particular, we would like to
define gauge transformations in such a way that a gauge
transformation of a connection one-form is still a connection
one-form. A couple of remarks on this problem may be found in
[3]. Also, it would be interesting to equip
our algebraic constructions with a some kind of topology,
like $C^ *$ or Frechet topology. Some topological aspects
of quantum fibre bundles are discussed in [8].
Furthermore, the theory of quantum
fibre bundles reviewed in this article is strictly related to
the theory of algebraic extensions. We think that the analysis
of quantum bundles from the point of view of Hopf-Galois
extensions may lead to a deeper insight into the both subjects.
Finally, we think it is desirable to develop generalised
fibre bundles defined in Remark~2.2.5. in order to construct
a gauge theory on general homogeneous spaces. The developement
of such a theory becomes even more important and challenging now
that the appearence of the $SU_q(2)$ homogeneous spaces
in the Connes description of Standard Model was announced [10].

\section{Acknowledgements}
Most of the results presented in this paper were
obtained jointly with Shahn Majid.  I would like to thank him
for a fruitful collaboration and many interesting
discussions. This paper was written during my
stay at the Universite Libre de Bruxelles;  I am grateful to the European Union
for the fellowship in the framework
of the Human Capital and Mobility Scheme. My work was also
supported by the grant KBN 2 P 302 21706 p 01.
\references{No}{\item{[1]}
R.J. \spa{Blattner,}  M. \spa{Cohen} and S. \spa{Montgomery},
{\it Crossed Products and
Inner Actions of Hopf Algebras}\/,  Trans. Amer.
Math. Soc.   298 (1986) 671.
\item{[2]} T. \spa{Brzezi{\'n}ski}, {\it
Differential
Geometry of Quantum Groups and Quantum Fibre Bundles}\/,
University of Cambridge, PhD thesis, 1994.
\item{[3]} T. \spa{Brzezi{\'n}ski}, {\it Remarks on Quantum Principal
Bundles}\/, In. {\it Quantum Groups. Formalism and
Applications}, J. Lukierski, Z. Popowicz and J. Sobczyk, eds.
Polish Scientific Publishers PWN, 1995, p. 3.
\item{[4]} T. \spa{Brzezi{\'n}ski}, {\it Translation Map in Quantum
Principal Bundles}\/, preprint (1994) {\tt hep-th/9407145}.
\item{[5]} T. \spa{Brzezi{\'n}ski}, {\it Quantum Homogeneous
Spaces as Quantum Quotient Spaces}\/, in preparation (1995)
\item{[6]} T. \spa{Brzezi{\'n}ski}
and S. \spa{Majid}, {\it Quantum Group Gauge Theory on Quantum Spaces}\/,
Commun.~Math.~Phys.   157 (1993) 591; {it ibid.}
 167 (1995) 235 (erratum).
\item{[7]} T. \spa{Brzezi{\'n}ski}
and S. \spa{Majid}, {\it Quantum Group Gauge Theory on Classical
 Spaces}\/,
Phys. Lett.   B298 (1993) 339.
\item{[8]} R.J. \spa{Budzy{\'n}ski}  and
W. \spa{Kondracki}, {\it
Quantum principal fiber bundles: topological aspects}\/,
preprint (1994) {\tt hep-th/9401019}.
\item{[9]} A. \spa{Connes}, {\it Non-Commutative Geometry}\/,
Academic Press, 1994.
\item{[10]} A. \spa{Connes}, {\it A lecture given at the
Conference on Non-commutative Geometry and Its Applications}\/,
Castle T\v re\v s\v t, Czech Republic, May 1995.
\item{[11]} A. \spa{Connes} and M. \spa{Rieffel}, {\it
 Yang-Mills for Non-Commutative
Two-Tori}\/,  Contemp. Math. 62 (1987) 237.
\item{[12]}C.-S. \spa{Chu},
P.-M. \spa{Ho} and H. \spa{Steinacker}, {\it
Q-deformed {D}irac monopole
with arbitrary charge}\/, preprint (1994) {\tt hep-th/9404023}.
\item{[13]} Y. \spa{Doi}, {\it Equivalent Crossed Products
for a Hopf Algebra}\/,
Commun. Algebra 17 (1989) 3053.
\item{[14]} Y. \spa{Doi} and M. \spa{Takeuchi}, {\it
 Cleft Module Algebras
and Hopf Modules}\/,  Commun. Algebra 14 (1986) 801.
\item{[15]} V.G. \spa{Drinfeld}, {\it Quantum Groups}\/, In {\it
Proceedings of the
International Congress of
Mathematicians, Berkeley, Cal. Vol.1}, Academic Press, 1986, p.798.
\item{[16]} T. \spa{Eguchi,} P. \spa{Gilkey} and A. \spa{Hanson},
{\it  Gravitation,
Gauge Thoeries and Differential Geometry}\/, Phys. Rep.
 66 (1980) 213.
\item{[17]} P.M. \spa{Hajac}, {\it Strong Connections
and $U\sb
q(2)$-Yang-Mills Theory on Quantum Principal Bundles}\/,
preprint(1994) {\tt hep-th/9406129}.
\item{[18]} D. \spa{Husemoller}, {\it Fibre
Bundles}\/,
Springer-Verlag, 3rd ed. 1994.
\item{[19]} D. \spa{Kastler}, {\it Cyclic Cohomology
within Differential Envelope}\/, Hermann, 1988.
\item{[20]}E. \spa{Kunz},
{\it K\"ahler Differentials}\/, Vieweg \& Sohn, 1986.
\item{[21]} S. \spa{Majid}, {\it Cross Product Quantisation,
 Nonabelian
Cohomology and Twisting of Hopf Algebras}\/,
In. {\it Generalised Symmetries in Physics}, H.-D. Doebner,
V.K. Dobrev and A.G. Ushveridze, eds., World Scientific, 1994, p. 13.
\item{[22]} U. \spa{Meyer}, {\it Projective Quantum Spaces}\/,
 Lett. Math.
Phys. to appear.
\item{[23]}
M. \spa{Pflaum}, {\it Quantum Groups on Fibre Bundles}\/,
 Commun. Math.
Phys. 166 (1994) 279.
\item{[24]} P. \spa{Podle{\'s}}, {\it Quantum Spheres}\/, Lett. Math.
Phys. 14 (1087) 193.
\item{[25]} H.-J. \spa{Schneider} {\it Principal Homogeneous Spaces
for Arbitrary Hopf Algebras}\/, Israel J. Math. 72
(1990) 167; R.J. \spa{Blattner} and S. \spa{Montgomery}, {\it
 Crossed Products and Galois
Extensions of Hopf Algebras}\/, Pacific J. Math. 137
(1989) 37.
\item{[26]} M.E. \spa{Sweedler}, {\it Hopf Algebras}\/,
Benjamin, 1969.
\item{[27]} M.E. \spa{Sweedler}, {\it Cohomology of Algebras
over Hopf
Algebras}\/,  Trans. AMS 133 (1968) 205.
\item{[28]} S.L. \spa{Woronowicz}, {\it Twisted $SU\sb 2$ Group. An
Example of a Non-commutative Differential Calculus}\/,
Publ. RIMS Kyoto University 23 (1987) 117.
\item{[29]} S.L. \spa{Woronowicz}, {\it Differential Calculus on Compact
Matrix Pseudogroups (Quantum Groups)}\/,
Commun.~Math.~Phys.
 122 (1989) 125.
}
\bye